\documentclass[12pt]{article}
\usepackage[table]{xcolor}
\usepackage[T1]{fontenc}
\usepackage[utf8]{inputenc}

\usepackage{amsmath,amssymb}
\usepackage{fancybox}
\usepackage{tabularx}
\usepackage[normalem]{ulem}
\usepackage[dvips]{graphicx}
\usepackage{textcomp} 

\usepackage{pst-plot,pst-tree,pst-node,pstricks-add}

\usepackage{nopageno}
\usepackage{url}

\usepackage{theorem}

\usepackage{vmargin}  
\setmarginsrb{2cm}{2cm}{2cm}{2cm}{0cm}{0cm}{0cm}{0cm}

\usepackage[frenchb]{babel}
\usepackage[np]{numprint}
\begin{document}
\setlength\parindent{0mm}



\title{Résolution du jeu de Juniper Green}
\author{Julien Lemoine}
\maketitle

\section*{Résumé}

Dans cet article, nous décrivons la résolution du jeu de Juniper Green pour tout entier, avec des techniques suffisamment élémentaires pour être expliquées à des élèves maîtrisant les notions de diviseurs et multiples. Une méthode générale permet de résoudre le jeu pour presque tout entier, si bien qu'il ne reste qu'une petite cinquantaine de cas à traiter, chacun pouvant être résolu au prix d'une étude sur papier de quelques minutes.

\section*{Introduction}

En 1997, dans sa rubrique de mathématiques récréatives de Scientific American, Ian Stewart présente le jeu de Juniper Green \cite{is97}, un jeu qui se joue avec tous les entiers entre 1 et un certain $n$. Il explique la stratégie gagnante lorsque $n=40$, et laisse au lecteur le soin de résoudre le jeu avec $n=100$. À la fin de son article, il pose finalement la question générale : peut-on déterminer la stratégie gagnante pour tout entier $n$ ?

Lorsque je m'attelle au problème, quelques recherches sont disponibles en ligne. On peut assez facilement obtenir une stratégie qui résout presque tout entier $n$, la plus grande valeur de $n$ résistant à cette approche étant 118. De l'autre côté du spectre, quelques essais de programmation calculent les premières valeurs de $n$ (par exemple, Laval et Sicard donnent les résultats jusqu'à $n=50$ \cite{ls17}). Il reste à faire le lien entre 50 et 118, ce qui n'est a priori pas simple compte-tenu de l'explosion combinatoire.

Après avoir essayé de résoudre un maximum de valeurs avec une stratégie manuelle (c'est le cas des lignes en vert dans la table \ref{table_res}), je me lance dans la programmation. En implémentant uniquement deux simplifications (la suppression des composantes connexes inutiles dans le graphe décrivant la position, et la suppression des sommets de degré 1 et de leur voisin), on obtient un programme qui résout presque toutes les valeurs restantes en un temps raisonnable. Seuls restent $n=110$, 116 et 117, que quelques optimisations moins évidentes me permettent de finalement calculer. Évidemment, dès le lendemain, je découvre une méthode élémentaire qui aurait pu m'économiser des dizaines d'heures de temps de calcul... mais je doute que j'aurais pu y penser, sans me confronter à la programmation et à l'observation des situations problématiques.

Cette méthode de résolution est explicable à un élève de collège, les quelques calculs qu'il faut mener sur certaines valeurs de $n$ (celles qui sont imperméables à la stratégie générale) sont faisables en quelques minutes pour chaque valeur étudiée, avec un papier et un crayon et quelques essais/erreurs, et s'avèrent un bon exercice pour travailler les diviseurs/multiples et le calcul mental.

\section{Règle du jeu}
\label{regle}

Le jeu de Juniper Green, popularisé par Ian Stewart \cite{is97}, et dont la paternité est attribuée à Rob Porteous  \cite{is972}, est un jeu combinatoire impartial, qui se joue avec les entiers de $[|1;n|]$. Alice choisit un entier, puis Bob choisit un diviseur ou un multiple de ce nombre, puis Alice un diviseur ou un multiple de ce dernier nombre... Un nombre donné ne peut être choisi qu'une seule fois, et celui qui ne peut plus jouer a perdu.

Voici un exemple de partie avec $n=8$ : Alice choisit 2, Bob choisit 8, Alice choisit 4, Bob choisit 1, Alice choisit 7. Bob ne pouvant plus jouer, c'est Alice qui a gagné.

Enfin, comme nous le verrons dans la section \ref{élémentaire}, on ajoute une règle qui énonce que le premier joueur doit choisir un nombre pair lors du premier coup de la partie, car sinon, il dispose trop facilement d'une certaine stratégie gagnante.

\section{Représentation des positions sous forme de graphe}

On appelle \emph{position} l'état de la partie à un instant donné ; lorsqu'un joueur joue un \emph{coup}, il réalise une transition entre deux positions.

\begin{figure}[ht]
\centering
\includegraphics[scale=0.5]{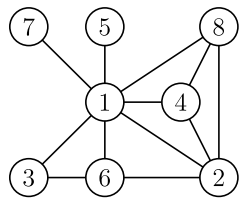}
\caption{Graphe initial du jeu avec $n=8$.}
\label{graphe_8}
\end{figure}

Il est possible de visualiser une position avec un graphe non orienté, dont les sommets sont les entiers de $[|1;n|]$ non encore choisis, et où deux sommets sont reliés par une arête si l'un est un multiple de l'autre. La figure \ref{graphe_8} présente ainsi le graphe initial du jeu avec $n=8$, c'est-à-dire l'état de la partie avant que le premier joueur ne joue.

\begin{figure}[ht]
\centering
\includegraphics[width=\textwidth]{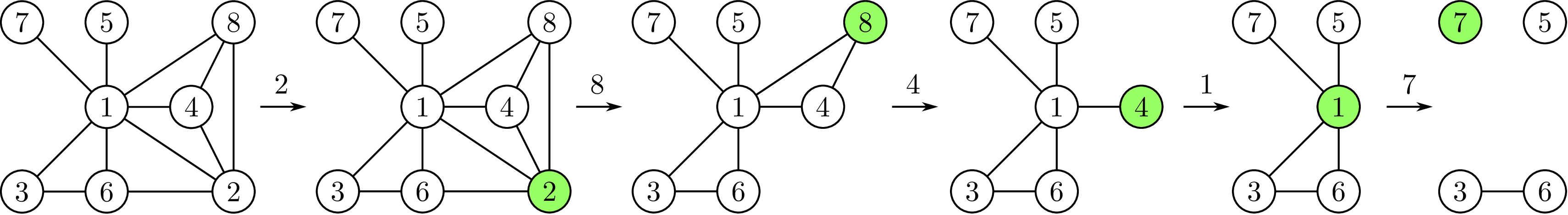}
\caption{Partie de la section \ref{regle}.}
\label{exemple_jeu}
\end{figure}

Le premier joueur choisit ensuite un nombre : on convient alors, dans la représentation, de colorier le sommet qui vient d'être choisi. Selon ce principe, la figure \ref{exemple_jeu} représente la partie de la section \ref{regle} ; les coups sont représentés par les flèches.

Remarquons qu'une fois le graphe construit, il est inutile de connaître le numéro des sommets. Le jeu de Juniper-Green se ramène à un simple jeu de déplacement sur un graphe non orienté, où on ne peut passer deux fois par le même sommet.

\section{Une première stratégie élémentaire}
\label{élémentaire}

Reprenant la figure \ref{graphe_8}, on voit facilement une stratégie victorieuse plus rapide pour Alice : elle choisit initialement 5 ; Bob est alors forcé de choisir 1, puis Alice choisit 7 et elle gagne. Cette stratégie est généralisable à tout $n\geq 3$ hormis 4, 6 et 10, c'est-à-dire dès lors qu'il y a au moins 2 nombres premiers dans $]\frac{n}{2};n]$\footnote{Démonstration en annexe \ref{resnbpr}.} (cela ne marcherait pas avec un nombre premier plus petit que $\frac{n}{2}$, car Bob pourrait choisir un multiple de ce nombre).

Pour cette raison, on impose au premier joueur la règle supplémentaire : le premier nombre qu'il choisit doit être pair. Ceci l'empêche de commencer par un nombre premier, et donc d'appliquer cette stratégie. Dans la suite, nous considérerons donc le jeu où cette règle est en vigueur.

\section{Simplification des positions}
\label{simplification}

Étant donnée une position du jeu, on dira que son \emph{résultat} est \emph{gagnante} si le joueur dont c'est le tour dispose d'une stratégie lui assurant la victoire, \emph{perdante} sinon.

Une position est gagnante si, à partir d'elle, on peut trouver un coup qui renvoie l'adversaire à une position perdante ; symétriquement, une position est perdante si, à partir d'elle, tous les coups renvoient l'adversaire à une position gagnante.

Par exemple, sur la figure \ref{exemple_jeu}, la position du milieu (avec le sommet 8 actif) est gagnante : si on joue 4, l'autre répond 1 (coup forcé), puis on joue 5 et on a gagné. L'existence d'un coup conduisant à une position perdante suffit à assurer le statut gagnant de la position (même si l'autre coup, 1, mène lui au contraire à une position gagnante pour l'adversaire, et doit donc être évité).

Inversement, la position à la gauche de celle-ci (avec le sommet 2 actif) est perdante :
\begin{itemize}
\item Si on joue 1, l'autre joue 5 et on a perdu ;
\item Si on joue 4, l'autre joue 8, le coup 1 est forcé, l'autre joue 5 et on a perdu ;
\item Symétriquement, si on joue 8, l'autre joue 4 et la partie se termine de même.
\end{itemize}
La position est perdante car tous les coups possibles mènent à une position gagnante.

\begin{figure}[ht]
\centering
\includegraphics[scale=0.5]{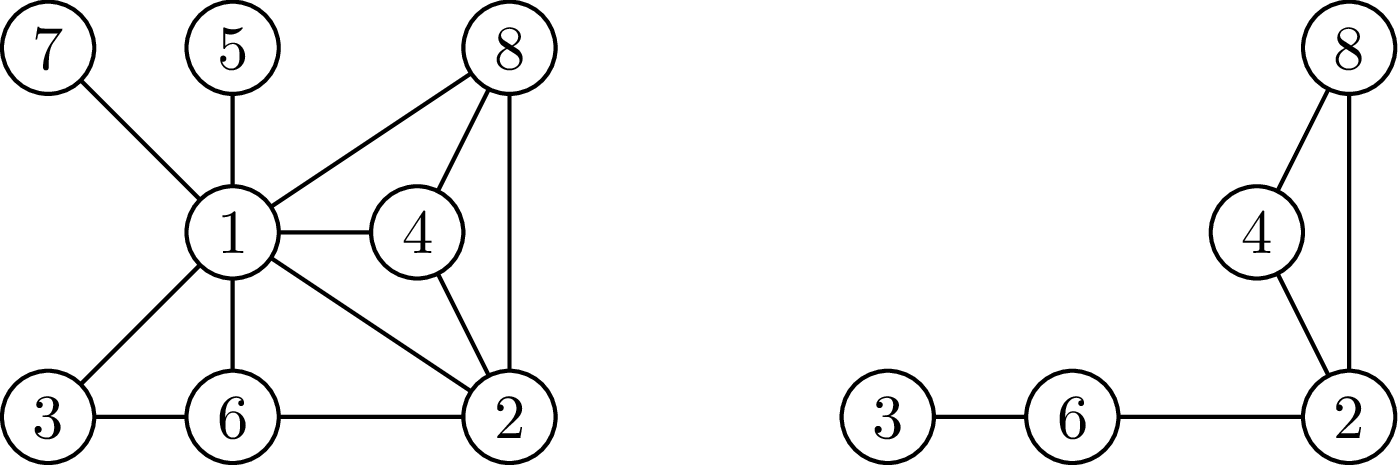}
\caption{Ces deux positions ont le même résultat.}
\label{graphe_8_simpl}
\end{figure}

On dira qu'on \emph{simplifie} une position si, en retirant des sommets, la position obtenue a le même résultat ; on peut en effet parfois prédire cela sans même être qu'il soit nécessaire de calculer ce résultat. Par exemple, si on reprend la position de la figure \ref{graphe_8}, et si on enlève les sommets 1, 5 et 7, on obtient une position de même résultat (voir figure \ref{graphe_8_simpl}). En effet, dans le jeu avec les sommets 1, 5, 7, personne n'a intérêt à jouer le sommet 1, car l'autre répondrait immédiatement 5 pour gagner la partie. Donc, on n'accède à ces sommets qu'en dernier recours. Ainsi, sur la figure \ref{graphe_8_simpl}, le joueur qui perd pour la position de droite perd aussi pour celle de gauche : il essaie autant que possible de ne pas jouer 1 en jouant comme sur la position de droite, et dès lors qu'il est bloqué avec les sommets de la position de droite, il ne lui reste plus qu'à jouer 1 et perdre au tour suivant.

En général, pour le jeu avec les $n$ premiers entiers, on simplifiera de manière analogue 1 et les nombres premiers strictement supérieurs à $\frac{n}{2}$. D'autres simplifications existent, mais celle-ci sera suffisante pour notre étude.

\section{Résolution manuelle pour presque tout $n$}
\label{presque}

Nous allons maintenant présenter des méthodes permettant de résoudre presque toutes les valeurs de $n$ sans l'aide de l'ordinateur. Ces méthodes sont issues d'idées présentées dans la résolution de $n=40$ décrite dans l'article de Stewart\cite{is97}, idées qui furent ensuite approfondies par Daniel Djament dans le bulletin de l'APMEP \cite{dd04} ou par divers contributeurs du site \url{diophante.fr} \cite{dio442} \cite{dio453}.

\begin{figure}[ht]
\centering
\includegraphics[scale=0.5]{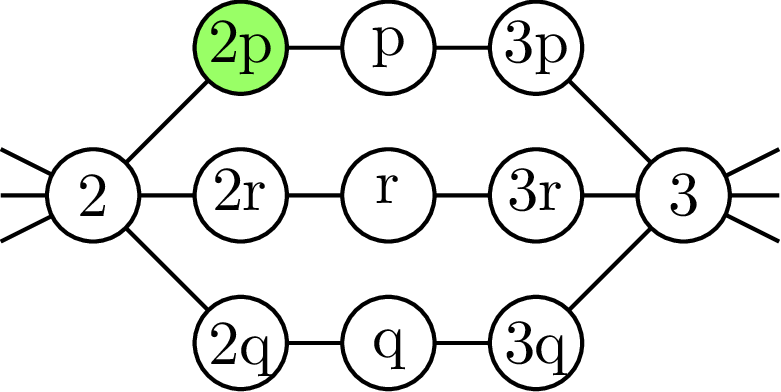}
\caption{Résolution de certaines valeurs de $n$.}
\label{graphe_pqr}
\end{figure}

Nous allons commencer par décrire une stratégie qui fonctionne dès lors qu'il existe au moins 3 nombres premiers dans $]\frac{n}{4};\frac{n}{3}]$ ; on en sélectionne alors trois quelconques que l'on nomme $p$, $q$ et $r$. On a tracé sur la figure \ref{graphe_pqr} une partie du graphe initial. La condition $p\in]\frac{n}{4};\frac{n}{3}]$ fait que les sommets $p$, $2p$ et $3p$ sont de degré 2 (et de même pour les sommets $q$ et $r$), après simplification du sommet 1.

Décrivons maintenant une stratégie gagnante pour le premier joueur. Il commence par choisir le sommet $2p$, comme illustré sur la figure \ref{graphe_pqr}, puis lorsqu'il atteint l'un des deux sommets 2 ou 3, il dirige son adversaire vers l'autre de ces deux sommets en passant par le sommet $q$. Enfin, une fois l'autre de ces deux sommets atteint, il dirige son adversaire vers le sommet $r$. Ainsi, son adversaire se retrouve bloqué après le parcours des 3 multiples de $r$.

\begin{figure}[ht]
\centering
\includegraphics[scale=0.5]{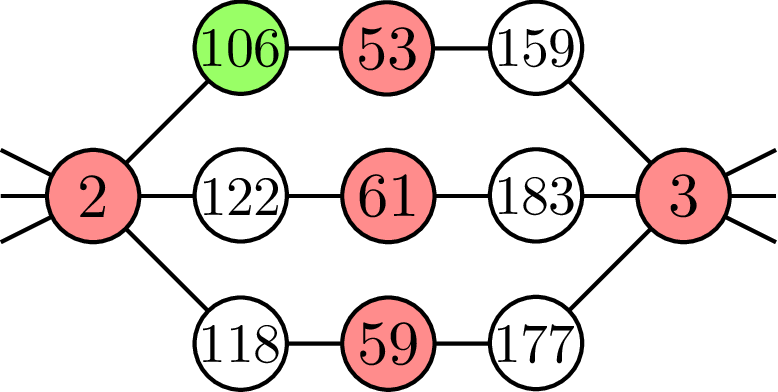}
\caption{Résolution de $n=200$.}
\label{strategie200}
\end{figure}

La figure \ref{strategie200} illustre ainsi la résolution de $n=200$. Le premier joueur commence avec 106. Lorsque c'est le tour de son adversaire, celui-ci est systématiquement placé sur un sommet blanc, qui étant de degré deux, ne lui laisse aucun choix, sauf lors de son premier coup, où il a le choix entre 2 et 53. Quel que soit ce choix, le premier joueur lui fait alors faire le tour du cycle externe, et une fois que 2 et 3 ont été joués, il le dirige vers le chemin du milieu, ce qui lui permet de finir la partie. La partie est donc soit 106, 2, 118, 59, 177, 3, 183, 61, 122, soit 106, 53, 159, 3, 177, 59, 118, 2, 122, 61, 183, si bien que tôt ou tard, le deuxième joueur sera condamné à jouer 1 puis à perdre au coup suivant.

La condition de l'existence de 3 nombres premiers dans $]\frac{n}{4};\frac{n}{3}]$ est vérifiée pour tout $n$ sauf $[|1;110|]$, $[|116;122|]$, $[|124;128|]$ et $[|172;176|]$\footnote{Démonstration en annexe \ref{resnbpr}.}. Ce sont donc les seules valeurs de $n$ qu'il nous reste à résoudre, par d'autres moyens.

\section{Autre stratégie similaire}
\label{similaires}

Une autre stratégie du même type que celle de la section \ref{presque} permet de résoudre certaines des valeurs de $n$ restantes en montrant là aussi que le premier joueur est en position de gagner. Elle nécessite qu'il y ait 2 nombres premiers dans $]\frac{n}{4};\frac{n}{3}]$, cette condition n'étant pas suffisante.

\begin{figure}[ht]
\centering
\includegraphics[scale=0.5]{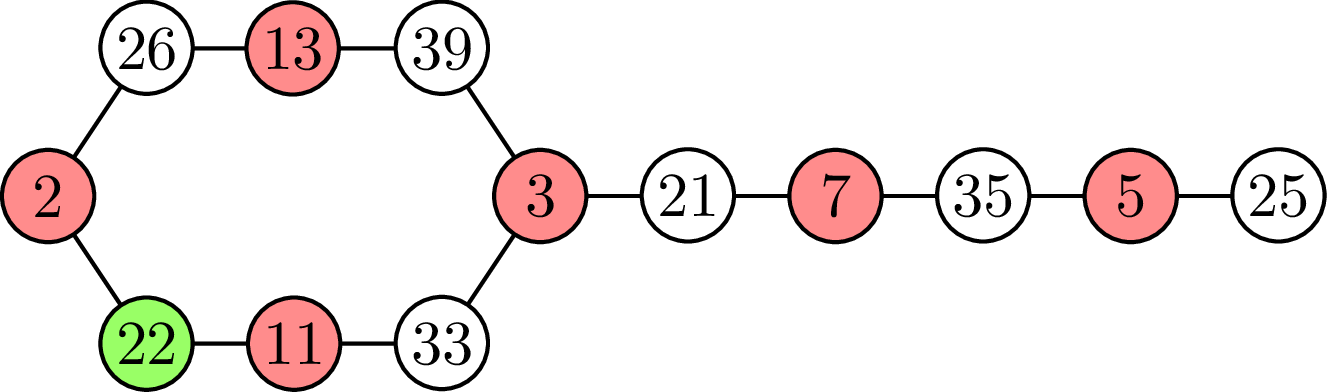}
\caption{Stratégie résolvant $n\in[|39, 41|]$.}
\label{strategie40}
\end{figure}

Cette stratégie repose sur le principe observable sur la figure \ref{strategie40} : les deux nombres premiers $p$ et $q$ de $]\frac{n}{4};\frac{n}{3}]$ (ici, 11 et 13) créent un cycle avec 2 et 3. Le premier joueur choisit d'abord $2p$ (ou $2q$), soit ici 22 (ou 26). Dans la suite, il joue toujours un produit de deux nombres premiers. Les nombres qu'il joue doivent être assez grands pour qu'ils n'aient pas de multiple inférieur à $n$ (ici, ils doivent valoir au moins 21), ainsi le deuxième joueur est forcé de répondre en jouant des nombres premiers.

Hormis le premier coup où il dispose de deux choix, les coups du deuxième joueur sont des coups forcés, et le cycle à gauche de la figure fait qu'il va devoir jouer 3. Une fois le coup 3 joué par le deuxième joueur, le premier joueur joue 21, et la suite de coups forcés se termine par 25 qui scelle la victoire du premier joueur. Ici, le coup terminal est 25, c'est-à-dire $5^2$ ; en général, c'est un coup du type $rs$, où $r$ et $s$ sont deux nombres premiers qui ont déjà été joués (avec éventuellement comme ici, $r=s$).

Cette stratégie particulière ne fonctionne pas pour $n\leq 38$, car 39 apparaît dans l'étude ; elle ne fonctionne pas non plus pour $n \geq 42$, car le deuxième joueur dispose alors d'une échappatoire lorsque le premier joueur joue 21. En généralisant, elle est valable pour $n\in[|a, 2b-1|]$, où $a$ est le plus grand nombre utilisé par le premier joueur, et $b$ le plus petit. Notant $a$ et $b$ en gras, toute cette étude se résume avec la notation suivante :
$$n\in[|39, 41|] : [2~22~11~33~3][2~26~13~\textbf{39}~3][3~\textbf{21}~7~35~5~25]$$
Avec la même idée, nous pouvons alors résoudre d'autres valeurs de $n$ :

\begin{itemize}
\item $n\in[|57, 65|] : [2~34~17~51~3][2~38~19~\textbf{57}~3][3~\textbf{33}~11~55~5~35~7~49]$
\item $n\in[|65, 67|] : [2~\textbf{34}~17~51~3][2~38~19~57~3][3~39~13~\textbf{65}~5~35~7~49]$
\item $n\in\{69\} : [2~38~19~57~3][2~46~23~\textbf{69}~3][3~39~13~65~5~\textbf{35}~7~49]$
\item $n\in[|87, 91|] : [2~\textbf{46}~23~69~3][2~58~29~\textbf{87}~3][3~51~17~85~5~55~11~77~7~49]$
\item $n\in[|93, 101|] : [2~58~29~87~3][2~62~31~\textbf{93}~3][3~\textbf{51}~17~85~5~55~11~77~7~91~13~65]$
\item $n\in[|95, 109|] : [2~58~29~87~3][2~62~31~93~3][3~57~19~\textbf{95}~5~\textbf{55}~11~77~7~91~13~65]$
\item $n\in[|119, 123|] : [2~\textbf{62}~31~93~3][2~74~37~111~3][3~69~23~115~5~65~13~91~7~\textbf{119}~17~85]$
\item $n\in[|123, 129|] : [2~74~37~111~3][2~82~41~\textbf{123}~3][3~69~23~115~5~\textbf{65}~13~91~7~119~17~85]$
\item $n\in[|161, 185|] : [2~94~47~141~3][2~106~53~159~3][3~\textbf{93}~31~155~5~115~23~\textbf{161}~7~133~19~95]$
\end{itemize}

\medskip

Faisons un bilan des valeurs de $n$ résolues jusqu'ici. Ces résolutions manuelles apparaissent en vert dans la table \ref{table_res}, et il ne reste désormais à résoudre que les valeurs de $n$ suivantes : 
$$[|1;39|]~~[|42;56|]~~68~~[|70;86|]~~92~~110~~[|116;118|]$$

Il ne sera toutefois pas nécessaire de réaliser un calcul pour l'ensemble de ces valeurs, car nous allons montrer que plusieurs de ces calculs sont équivalents.

\section{Valeurs de $n$ avec le même résultat}

Deux remarques permettent de limiter le nombre de calculs à mener parmi les valeurs de $n$ qu'il nous reste à étudier. Tout d'abord, si $p\geq 5$ est un nombre premier, alors $n=p-1$ et $n=p$ ont le même résultat. En effet, il y a alors au moins un autre nombre premier $q > p/2$, jouable uniquement à partir de 1. Jouer $p$ ou jouer $q$ revient au même et finit la partie, ajouter la possibilité de jouer $p$ ne change donc rien.

Par ailleurs, si $p\geq 3$ est un nombre premier, alors $n=2p-1$ et $n=2p$ ont le même résultat. En effet, un joueur qui jouerait $2p$ se verrait répondre $p$, il devrait alors jouer $1$ et il perdrait au coup suivant. Ajouter la possibilité de jouer $2p$ ne change rien là non plus.

Nous avons indiqué par le symbole « | » les valeurs de $n$ pour lesquelles ces remarques permettent de déduire qu'elles ont le même résultat que $n-1$. Parfois, plus de deux valeurs de $n$ successives sont associées : par exemple, la première remarque associe $n=4$ et $n=5$, mais aussi $n=6$ et $n=7$, tandis que la deuxième associe $n=5$ et $n=6$, d'où on déduit que ces 4 valeurs de $n$ ont le même résultat ; il suffit donc de mener le calcul pour $n=4$.

\section{Stratégie d'appariement}
\label{paires}

Une idée simple permet de résoudre les valeurs de $n$ restantes : une position est perdante dès lors que l'on peut partitionner les entiers de la position en paires de sommets adjacents. En effet, dès lors que le premier joueur joue un entier, le deuxième joueur joue l'autre entier de la paire, et ce jusqu'à la fin de la partie.

\begin{figure}[ht]
\centering
\includegraphics[scale=0.5]{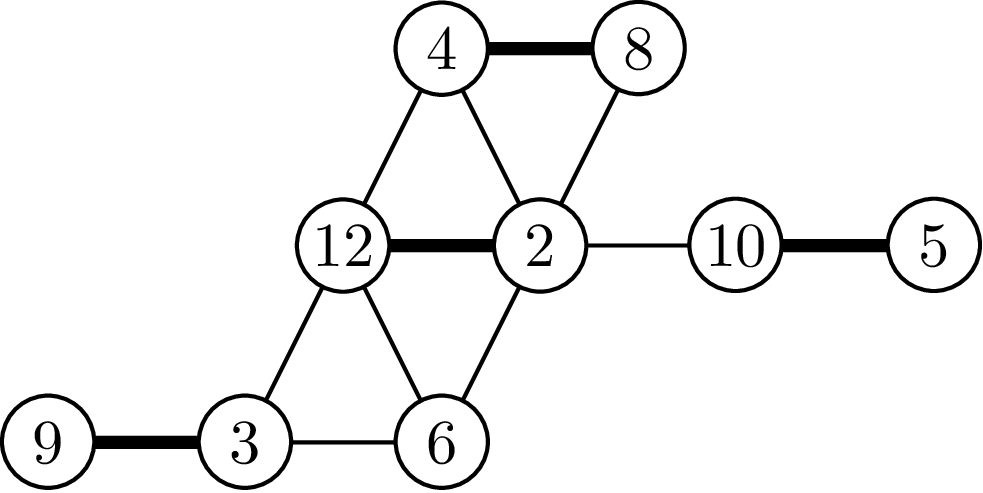}
\caption{Apppariements pour le jeu avec $n=12$.}
\label{graphe_12}
\end{figure}

Donnons l'exemple de l'étude du jeu avec $n=9$. On commence par éliminer 1, ainsi que 5 et 7 (les nombres premiers supérieurs strictement à $\frac{n}{2}$). Puis on associe les entiers restants par paires de diviseurs et multiples (en gras sur la figure \ref{graphe_12}) : 9 étant de degré 1, la paire $9-3$ est forcée, puis $6-2$ est également forcée, et il ne reste que $4-8$. On en déduit que la position de départ est perdante (noté «P» dans la table \ref{table_res}) : le premier joueur va jouer un certain nombre d'entiers, le deuxième joueur répondra systématiquement avec l'autre entier de chaque paire. Finalement, le premier joueur sera forcé de jouer 1, et le deuxième gagnera en répondant 5 (ou 7).

\begin{figure}[ht]
\centering
\includegraphics[scale=0.5]{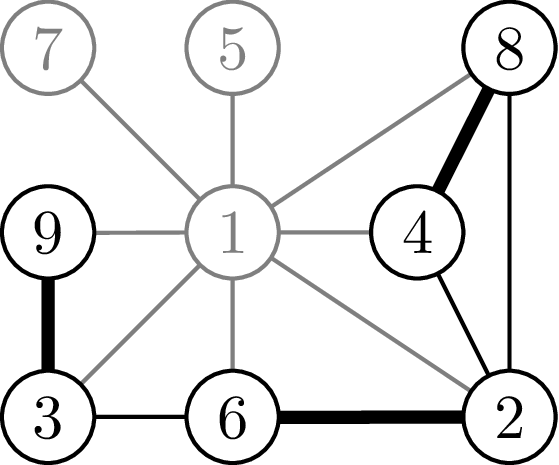}
\caption{Apppariements pour le jeu avec $n=9$.}
\label{graphe_9}
\end{figure}

Quand la position de départ est gagnante, un entier reste tout seul après les appariements, et c'est celui que joue le premier joueur en début de partie. Par exemple, pour $n=12$, on élimine 1, 7 et 11. Puis les paires $5-10$ et $9-3$ sont forcées, et enfin il reste 5 sommets. On peut par exemple apparier $4-8$ et $12-2$, de sorte que le premier joueur peut gagner en commençant par 6, puis en appliquant la même stratégie que précédemment.

Quand on étudie les autres valeurs de $n$ restantes, on voit que l'on dispose d'une certaine latitude pour déterminer les appariements : certains sont forcés, mais les sommets qui restent ensuite ont souvent un degré assez important, ce qui fait qu'il existe de nombreuses façons de les apparier, or trouver un seul appariement suffit.

De fait, la détermination de l'appariement est possible pour toutes ces valeurs de $n$, à l'exception notable de 49 et de 92 qui seront traités un peu plus loin. Une étude de quelques minutes suffit pour chaque valeur de $n$, même les plus grandes, pour peu qu'on s'y prenne avec méthode. Expliquons comment faire avec le cas de $n=56$.

Comme précédemment, on commence par éliminer 1, ainsi que 29, 31, 37, 41, 43, 47, 53. Puis, on procède aux appariements forcés : après l'élimination de 1, le sommet 49 n'est plus adjacent qu'à 7, d'où la paire $7-49$. De même, il y a les paires $23-46$, et $19-38$. Puis, du fait qu'on a déjà apparié 7, le sommet 35 n'est plus relié qu'à 5, d'où la paire $35-5$. Chaque paire déterminée est ainsi susceptible de provoquer un ou plusieurs autres appariements. On déduit ensuite $55-11$, $33-3$, $39-13$, $51-17$, $34-2$, $22-44$, $26-52$, $25-50$, et $21-42$.

Il reste alors les entiers 4, 6, 8, 9, 10, 12, 14, 15, 16, 18, 20, 24, 27, 28, 30, 32, 36, 40, 45, 48, 54, 56. On peut alors faire fonctionner un algorithme glouton, qui donne la priorité aux sommets de faible degré. Ainsi, parmi les entiers restants, le degré minimal est 2. On choisit un sommet de degré 2, par exemple 14, et l'un de ses voisins, par exemple 56. On a donc la paire $14-56$, et on continue ainsi : il vient les paires $28-4$, $15-30$, $45-9$, $27-54$, $10-40$, $18-6$, $36-12$, $24-8$ et $32-16$.

\begin{figure}[ht]
\centering
\includegraphics[scale=0.3]{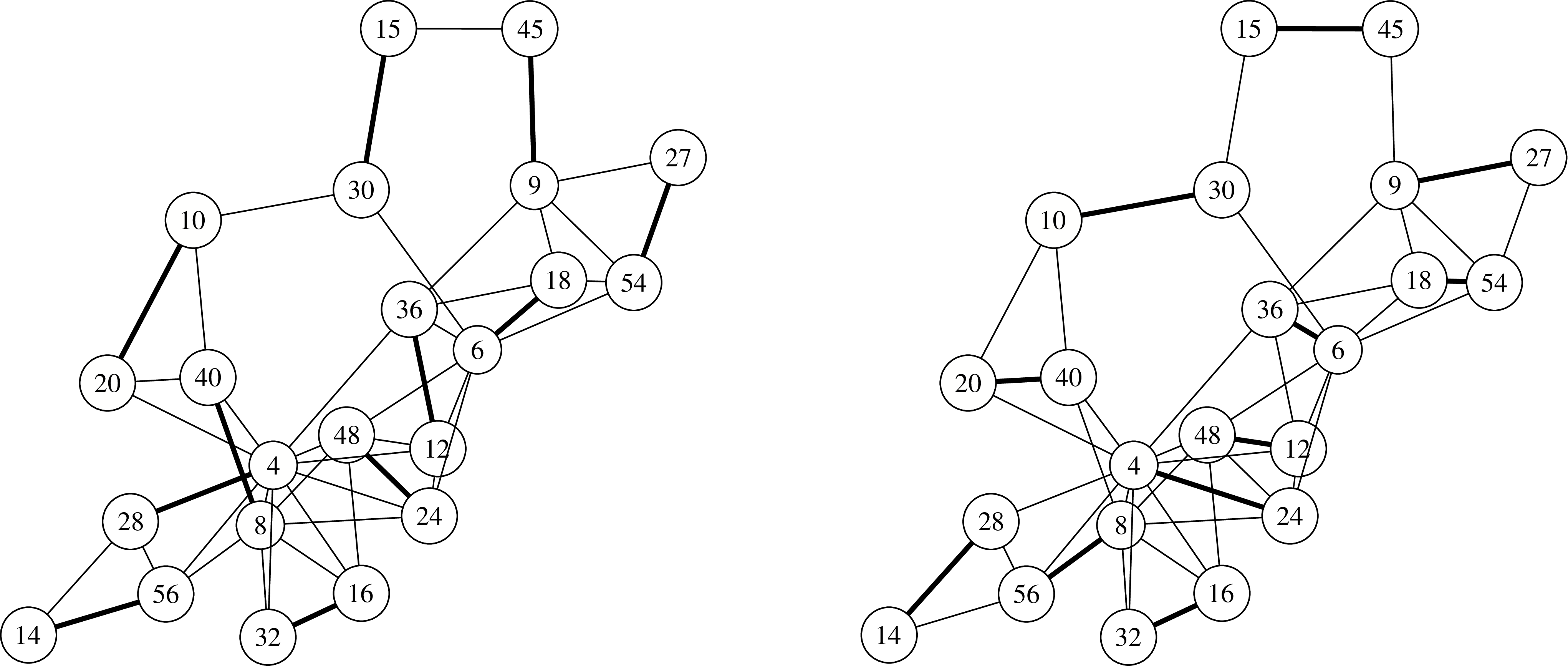}
\caption{Appariements pour le jeu avec $n=56$.}
\label{graphe_56}
\end{figure}

Mais il reste les entiers 20 et 48, qui ne sont pas multiples l'un de l'autre : il faut donc réorganiser les paires pour permettre à ces sommets de s'apparier. Ici, on peut remarquer qu'en cassant $10-40$ et $24-8$, on peut apparier les entiers en $10-20$, $8-40$ et $24-48$. Les autres valeurs de $n$ se traitent de manière similaire : s'il reste certains entiers non appariés à la fin de la procédure, il n'y a pas besoin de casser beaucoup de paires pour incorporer les entiers restants. Rappelons qu'il peut y avoir de multiples appariements possibles : l'appariement déterminé ici est présenté à gauche de la figure \ref{graphe_56}, un autre appariement qui convient est visible à droite de la figure.

\begin{figure}[ht]
\centering
\includegraphics[scale=0.5]{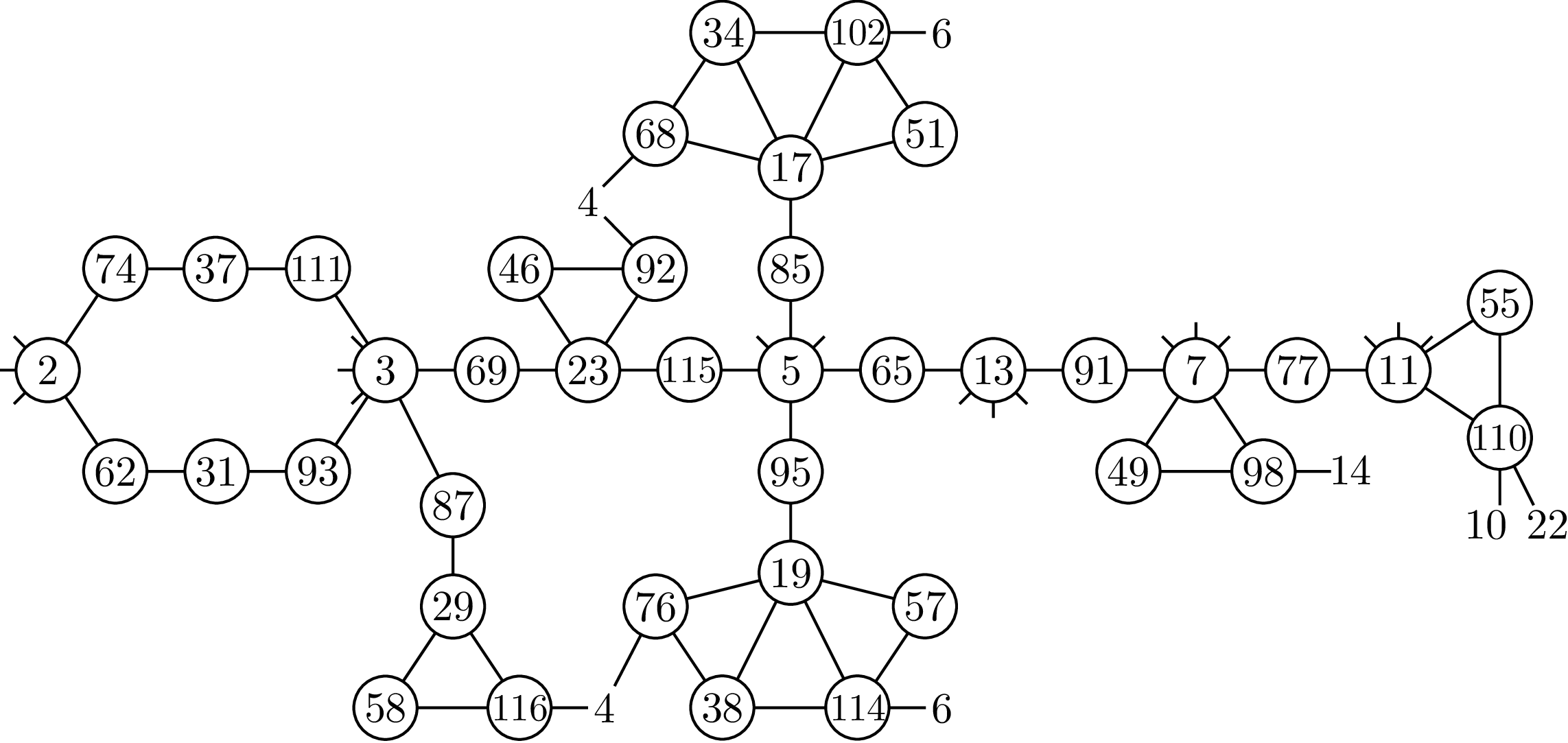}
\caption{Sommets s'appariant rapidement dans l'étude de $n=117$.}
\label{117}
\end{figure}

Nous donnons enfin un dernier exemple d'étude avec $n=117$. Après suppression de 1, des nombres premiers supérieurs à $\frac{n}{2}$, et des 4 paires forcées $53-106$, $47-94$, $43-86$ et $41-82$, il reste, entre autres, les sommets de la figure \ref{117}.

Le cycle à gauche de la figure doit alors engendrer 4 paires : ou bien $2-74$, $37-111$, $3-93$ et $31-62$, ou bien $2-62$, $31-93$, $3-111$ et $37-74$. Il n'est ensuite plus nécessaire de considérer 2 et 3, déjà appariés, c'est pourquoi les arêtes qui leur sont incidentes n'apparaissent pas sur le reste du graphe (comme avec 46, 51, etc).

L'appariement des autres sommets entourés de la figure est forcé (nous laissons au lecteur le soin de mettre en gras les arêtes sur la figure), et il reste alors 54 entiers à apparier, ce qui ne s'avère  pas très compliqué vu leurs très grands nombres de diviseurs et multiples.

La mise en œuvre de cette stratégie nous a permis de déterminer les appariements décrits dans l'annexe \ref{appa}, qui résolvent les valeurs de $n$ qui restaient en suspens.

\section{Cas particuliers : $n=49$ et $n=92$}

Les positions de départ se scindent en deux ensembles : soit le premier joueur force la victoire avec des sommets de faible degré, grâce à la méthode des parties \ref{presque} et \ref{similaires} ; soit les nombreuses arêtes du graphe font qu'on peut utiliser la méthode de l'appariement de la partie \ref{paires}.

Dans ce deuxième cas, on peut facilement prévoir le vainqueur : partant de $n$ sommets, on enlève 1, ainsi que les nombres premiers strictement supérieurs à $\frac{n}{2}$. Si le nombre de sommets restants est pair, on peut prédire que la position de départ sera perdante, sinon elle sera gagnante. Par exemple, pour le cas $n=56$, il y a 7 nombres premiers simplifiés initialement ; $56-1-7 = 48$, qui est pair, et la position de départ est perdante, car on peut apparier ces 48 sommets.

Cette règle permet de prévoir le vainqueur pour tout $n\geq 4$ qui n'est pas colorié en vert dans la table \ref{table_res}, car on peut toujours apparier les sommets qui restent après la simplification initiale, sauf dans deux cas : $n=49$ et $n=92$. En effet, dans ces deux cas, le joueur en position de force peut s'assurer qu'un certain sommet ne sera jamais joué, ce qui change la parité du nombre de coups joués, et modifie le résultat prévu.

Plus précisément, quand on résout $n=49$, après la simplification initiale, les appariements $49-7$ et $25-5$ sont forcés (car les sommets 49 et 25 sont alors de degré 1). Le sommet 35 se retrouve isolé et ne sera jamais joué. De même, pour l'étude de $n=92$, les appariements $49-7$, puis $77-11$, $55-5$ et $65-13$ sont forcés, si bien que $91=7\times 13$ ne sera jamais joué.

Des appariements qui conviennent pour $n=49$ et $n=92$ sont décrits dans l'annexe \ref{appa}.

\section{Résultats}

Les résultats obtenus sont récapitulés dans la table \ref{table_res} et sont cohérents avec les résultats obtenus jusqu'à $n=50$ par Laval et Sicard dans leur étude \cite{ls17}. Un « G » signifie que le premier joueur dispose d'une stratégie gagnante, un « P » qu'il commence depuis une position perdante.

\begin{table}[ht]
\centering
\begin{minipage}{0.16\textwidth}
\begin{tabular}{ |c|c| }
\hline
$n$ & G/P\tabularnewline
\hline
1 & -\tabularnewline
2 & P\tabularnewline
3 & G\tabularnewline
4 & P\tabularnewline
5 & |\tabularnewline
6 & |\tabularnewline
7 & |\tabularnewline
8 & G\tabularnewline
9 & P\tabularnewline
10 & |\tabularnewline
\hline
\end{tabular}
\end{minipage}
\begin{minipage}{0.16\textwidth}
\begin{tabular}{ |c|c| }
\hline
$n$ & G/P\tabularnewline
\hline
11 & |\tabularnewline
12 & G\tabularnewline
13 & |\tabularnewline
14 & |\tabularnewline
15 & P\tabularnewline
16 & G\tabularnewline
17 & |\tabularnewline
18 & P\tabularnewline
19 & |\tabularnewline
20 & G\tabularnewline
\hline
\end{tabular}
\end{minipage}
\begin{minipage}{0.16\textwidth}
\begin{tabular}{ |c|c| }
\hline
$n$ & G/P\tabularnewline
\hline
21 & P\tabularnewline
22 & |\tabularnewline
23 & |\tabularnewline
24 & G\tabularnewline
25 & P\tabularnewline
26 & |\tabularnewline
27 & G\tabularnewline
28 & P\tabularnewline
29 & |\tabularnewline
30 & G\tabularnewline
\hline
\end{tabular}
\end{minipage}
\begin{minipage}{0.16\textwidth}
\begin{tabular}{ |c|c| }
\hline
$n$ & G/P\tabularnewline
\hline
31 & |\tabularnewline
32 & P\tabularnewline
33 & G\tabularnewline
34 & |\tabularnewline
35 & P\tabularnewline
36 & G\tabularnewline
37 & |\tabularnewline
38 & |\tabularnewline
39 & \cellcolor{green!25}G\tabularnewline
40 & \cellcolor{green!25}G\tabularnewline
\hline
\end{tabular}
\end{minipage}
\begin{minipage}{0.16\textwidth}
\begin{tabular}{ |c|c| }
\hline
$n$ & G/P\tabularnewline
\hline
41 & \cellcolor{green!25}G\tabularnewline
42 & P\tabularnewline
43 & |\tabularnewline
44 & G\tabularnewline
45 & P\tabularnewline
46 & |\tabularnewline
47 & |\tabularnewline
48 & G\tabularnewline
49 & G\tabularnewline
50 & G\tabularnewline
\hline
\end{tabular}
\end{minipage}
\begin{minipage}{0.16\textwidth}
\begin{tabular}{ |c|c| }
\hline
$n$ & G/P\tabularnewline
\hline
51 & P\tabularnewline
52 & G\tabularnewline
53 & |\tabularnewline
54 & P\tabularnewline
55 & G\tabularnewline
56 & P\tabularnewline
57 & \cellcolor{green!25}G\tabularnewline
58 & \cellcolor{green!25}G\tabularnewline
59 & \cellcolor{green!25}G\tabularnewline
60 & \cellcolor{green!25}G\tabularnewline
\hline
\end{tabular}
\end{minipage}

\begin{minipage}{0.16\textwidth}
\begin{tabular}{ |c|c| }
\hline
$n$ & G/P\tabularnewline
\hline
61 & \cellcolor{green!25}G\tabularnewline
62 & \cellcolor{green!25}G\tabularnewline
63 & \cellcolor{green!25}G\tabularnewline
64 & \cellcolor{green!25}G\tabularnewline
65 & \cellcolor{green!25}G\tabularnewline
66 & \cellcolor{green!25}G\tabularnewline
67 & \cellcolor{green!25}G\tabularnewline
68 & G\tabularnewline
69 & \cellcolor{green!25}G\tabularnewline
70 & G\tabularnewline
\hline
\end{tabular}
\end{minipage}
\begin{minipage}{0.16\textwidth}
\begin{tabular}{ |c|c| }
\hline
$n$ & G/P\tabularnewline
\hline
71 & |\tabularnewline
72 & P\tabularnewline
73 & |\tabularnewline
74 & |\tabularnewline
75 & G\tabularnewline
76 & P\tabularnewline
77 & G\tabularnewline
78 & P\tabularnewline
79 & |\tabularnewline
80 & G\tabularnewline
\hline
\end{tabular}
\end{minipage}
\begin{minipage}{0.16\textwidth}
\begin{tabular}{ |c|c| }
\hline
$n$ & G/P\tabularnewline
\hline
81 & P\tabularnewline
82 & |\tabularnewline
83 & |\tabularnewline
84 & G\tabularnewline
85 & P\tabularnewline
86 & |\tabularnewline
87 & \cellcolor{green!25}G\tabularnewline
88 & \cellcolor{green!25}G\tabularnewline
89 & \cellcolor{green!25}G\tabularnewline
90 & \cellcolor{green!25}G\tabularnewline
\hline
\end{tabular}
\end{minipage}
\begin{minipage}{0.16\textwidth}
\begin{tabular}{ |c|c| }
\hline
$n$ & G/P\tabularnewline
\hline
91 & \cellcolor{green!25}G\tabularnewline
92 & P\tabularnewline
93 & \cellcolor{green!25}G\tabularnewline
94 & \cellcolor{green!25}G\tabularnewline
95 & \cellcolor{green!25}G\tabularnewline
96 & \cellcolor{green!25}G\tabularnewline
97 & \cellcolor{green!25}G\tabularnewline
98 & \cellcolor{green!25}G\tabularnewline
99 & \cellcolor{green!25}G\tabularnewline
100 & \cellcolor{green!25}G\tabularnewline
\hline
\end{tabular}
\end{minipage}
\begin{minipage}{0.16\textwidth}
\begin{tabular}{ |c|c| }
\hline
$n$ & G/P\tabularnewline
\hline
101 & \cellcolor{green!25}G\tabularnewline
102 & \cellcolor{green!25}G\tabularnewline
103 & \cellcolor{green!25}G\tabularnewline
104 & \cellcolor{green!25}G\tabularnewline
105 & \cellcolor{green!25}G\tabularnewline
106 & \cellcolor{green!25}G\tabularnewline
107 & \cellcolor{green!25}G\tabularnewline
108 & \cellcolor{green!25}G\tabularnewline
109 & \cellcolor{green!25}G\tabularnewline
110 & P\tabularnewline
\hline
\end{tabular}
\end{minipage}
\begin{minipage}{0.16\textwidth}
\begin{tabular}{ |c|c| }
\hline
$n$ & G/P\tabularnewline
\hline
111 & \cellcolor{green!25}G\tabularnewline
112 & \cellcolor{green!25}G\tabularnewline
113 & \cellcolor{green!25}G\tabularnewline
114 & \cellcolor{green!25}G\tabularnewline
115 & \cellcolor{green!25}G\tabularnewline
116 & G\tabularnewline
117 & P\tabularnewline
118 & |\tabularnewline
119 & \cellcolor{green!25}G\tabularnewline
120+ & \cellcolor{green!25}G\tabularnewline
\hline
\end{tabular}
\end{minipage}
\caption{Résultats pour tout $n$.}
\label{table_res}
\end{table}

\section*{Conclusion}

Dans son article \cite{is97}, Ian Stewart demandait s'il était possible de résoudre le jeu de Juniper Green pour toute valeur de $n$. À l'aide des quelques méthodes élémentaires décrites dans cet article, c'est désormais le cas. Dans le feedback paru quelques mois plus tard \cite{is972}, il indiquait : « Paul J. Blatz of Van Nuys, Calif., recalled that the game was discussed in a number theory course given at Princeton University by Eugene P. Wigner in the late 1930s. A criterion for winning play in all cases was provided. The answer for JG-n depends on the oddness or evenness of the powers of primes that occur when n! is factorized ». Au vu de la résolution que nous présentons, ce critère semble étrange, et il serait intéressant de retrouver des traces de ce cours.

\appendix

\section{Résultats sur les nombres premiers}
\label{resnbpr}

Dans cet article, nous avons utilisé deux résultats concernant l'existence de nombres premiers dans certaines intervalles. Ils sont assez simples à démontrer en utilisant des théorèmes sur la répartition des nombres premiers.

Par exemple, dans sa thèse de doctorat \cite{pd98}, Pierre Dusart démontre que pour tout réel $x\geq 3275$, il existe un nombre premier $p$ dans l'intervalle $\displaystyle{\left]x ; x\left(1 + \frac{1}{2\ln^2 x}\right)\right]}$. On peut remplacer cet intervalle par $\displaystyle{\left]x ; x.\frac{132}{131}\right]}$ en remarquant que la fonction $\displaystyle{x\mapsto \left(1 + \frac{1}{2\ln^2 x}\right)}$ est décroissante, et que l'image de 3275 est plus petite que $\displaystyle{\frac{132}{131}}$.

On en déduit, en remplaçant $x$ par $\displaystyle{\frac{n}{2}}$, qu'il y a au moins un nombre premier dans $\displaystyle{\left]\frac{n}{2} ; \frac{n}{2}.\frac{132}{131}\right]}$, puis en remplaçant $x$ par $\displaystyle{\frac{n}{2}.\frac{132}{131}}$, qu'il y en a au moins un dans $\displaystyle{\left]\frac{n}{2}.\frac{132}{131} ; \frac{n}{2}.\left(\frac{132}{131}\right)^2\right]}$. Comme $\displaystyle{\frac{n}{2}.\left(\frac{132}{131}\right)^2 < n}$, il y a donc au moins deux nombres premiers dans $\displaystyle{\left]\frac{n}{2} ; n\right]}$ dès que $n\geq 6550$.

De même, comme $\displaystyle{\frac{n}{4}.\left(\frac{132}{131}\right)^3 < \frac{n}{3}}$, il y a au moins trois nombres premiers dans $\displaystyle{\left]\frac{n}{4} ; \frac{n}{3}\right]}$ dès que $n\geq 13100$.

Pour démontrer les résultats voulus, il ne reste plus qu'à utiliser un programme pour vérifier par le calcul l'existence des nombres premiers recherchés lorsque $n$ est inférieur respectivement à 6550 et 13100. Merci à Guillaume Aubian pour son aide concernant cette annexe.

\section{Appariements}
\label{appa}

Nous présentons ici des appariements qui permettent de résoudre certaines valeurs de $n$ (après simplification de 1 et des nombres premiers strictement supérieurs à $\frac{n}{2}$).

\footnotesize
\begin{verbatim}
4   (2,4)

8   (3,6)(2,8) 4

9   (3,9)(2,6)(4,8)

12  (5,10)(3,9)(2,12)(4,8) 6

15  (7,14)(3,9)(5,15)(2,10)(6,12)(4,8)

16  (7,14)(3,9)(5,15)(2,10)(6,12)(4,8) 16

18  (7,14)(5,15)(2,10)(8,16)(4,12)(3,9)(6,18)

20  (7,14)(3,9)(5,15)(2,10)(6,18)(4,12)(8,16) 20

21  (7,21)(2,14)(4,12)(3,9)(6,18)(8,16)(5,15)(10,20)

24  (11,22)(7,21)(2,14)(3,9)(5,15)(10,20)(6,18)(4,12)(8,16) 24

25  (11,22)(5,25)(3,15)(9,18)(7,21)(2,14)(10,20)(6,24)(4,12)(8,16)

27  (11,22)(13,26)(5,25)(3,15)(7,21)(2,14)(10,20)(9,27)(6,18)(4,12)(8,16) 24

28  (11,22)(13,26)(5,25)(3,15)(7,21)(9,27)(2,10)(14,28)(6,18)(4,20)(12,24)(8,16)

30  (11,22)(13,26)(5,25)(3,15)(7,21)(9,27)(2,10)(14,28)(6,18)(4,20)(12,24)(8,16) 30

32  (11,22)(13,26)(5,25)(3,15)(7,21)(9,27)(2,30)(14,28)(6,18)(4,16)(12,24)(10,20)(8,32)

33  (11,22)(13,26)(5,25)(3,15)(7,21)(9,27)(2,30)(14,28)(6,18)(4,16)(12,24)(10,20)(8,32) 33

35  (13,26)(17,34)(5,25)(7,35)(3,21)(15,30)(9,27)(11,33)(2,22)(10,20)(14,28)(6,18)(4,12)
    (8,24)(16,32)

36  (13,26)(17,34)(5,25)(7,35)(3,21)(15,30)(9,27)(11,33)(2,22)(10,20)(14,28)(6,18)(4,12)
    (8,24)(16,32) 36

42  (17,34)(19,38)(5,25)(7,35)(11,33)(2,22)(13,26)(3,39)(15,30)(21,42)(14,28)(9,27)(10,40)
    (4,20)(16,32)(6,18)(12,36)(8,24)

44  (17,34)(19,38)(5,25)(7,35)(11,33)(2,22)(13,26)(3,39)(15,30)(21,42)(14,28)(9,27)(10,40)
    (4,20)(16,32)(6,18)(12,36)(8,24) 44

45  (17,34)(19,38)(5,25)(7,35)(13,39)(2,26)(14,42)(3,21)(9,27)(4,28)(11,33)(22,44)(15,45)
    (6,18)(10,30)(20,40)(12,36)(8,24)(16,32)

48  (17,34)(19,38)(23,46)(5,25)(7,35)(13,39)(2,26)(14,42)(3,21)(9,27)(4,28)(11,33)(22,44)
    (15,45)(6,18)(10,30)(20,40)(12,36)(8,24)(16,32) 48

49  (17,34)(19,38)(23,46)(5,25)(7,49)(13,39)(2,26)(14,42)(3,21)(9,27)(4,28)(11,33)(22,44)
    (15,45)(6,18)(10,30)(20,40)(12,36)(8,24)(16,32) 35,48

50  (17,34)(19,38)(23,46)(7,49)(5,35)(25,50)(13,39)(2,26)(14,42)(3,21)(9,27)(4,28)(11,33)
    (22,44)(15,45)(6,18)(10,30)(20,40)(12,36)(8,24)(16,32) 48

51  (19,38)(23,46)(7,49)(5,35)(25,50)(13,39)(2,26)(17,34)(3,51)(21,42)(14,28)(9,27)(11,33)
    (22,44)(15,45)(6,18)(10,30)(4,20)(12,36)(8,40)(24,48)(16,32)

52  (19,38)(23,46)(7,49)(5,35)(25,50)(13,39)(2,26)(17,34)(3,51)(21,42)(14,28)(9,27)(11,33)
    (22,44)(15,45)(6,18)(10,30)(4,20)(12,36)(8,40)(24,48)(16,32) 52

54  (19,38)(23,46)(7,49)(5,35)(25,50)(17,51)(2,34)(14,42)(3,21)(4,28)(11,33)(22,44)(13,39)
    (26,52)(15,30)(9,45)(27,54)(10,20)(6,18)(12,36)(8,40)(16,32)(24,48)

55  (19,38)(23,46)(7,49)(5,35)(25,50)(11,55)(3,33)(21,42)(13,39)(17,51)(2,34)(14,28)(22,44)
    (26,52)(15,30)(9,45)(27,54)(10,40)(4,20)(6,18)(12,36)(8,24)(16,32) 48

56  (19,38)(23,46)(7,49)(5,35)(25,50)(11,55)(3,33)(21,42)(13,39)(17,51)(2,34)(22,44)(26,52)
    (14,56)(4,28)(15,30)(9,45)(27,54)(10,20)(6,18)(12,36)(8,40)(24,48)(16,32)

68  (23,46)(29,58)(31,62)(7,49)(5,35)(25,50)(11,55)(13,65)(3,39)(33,66)(17,51)(19,57)(2,38)
    (22,44)(26,52)(34,68)(21,42)(9,63)(27,54)(15,45)(14,56)(4,28)(6,18)(12,36)(8,24)(16,48)
    (32,64)(10,30)(20,40) 60

70  (29,58)(31,62)(7,49)(19,57)(2,38)(23,46)(3,69)(13,39)(26,52)(17,51)(34,68)(5,65)(25,50)
    (35,70)(11,55)(33,66)(22,44)(21,42)(9,63)(27,54)(15,45)(14,56)(4,28)(6,18)(12,36)(8,24)
    (16,48)(32,64)(10,30)(20,40) 60

72  (29,58)(31,62)(7,49)(19,57)(2,38)(23,46)(3,69)(13,39)(26,52)(17,51)(34,68)(5,65)(25,50)
    (35,70)(11,55)(33,66)(22,44)(21,42)(9,63)(27,54)(15,45)(14,56)(4,28)(6,18)(30,60)(10,20)
    (36,72)(12,24)(8,40)(16,48)(32,64)

75  (29,58)(31,62)(37,74)(7,49)(19,57)(2,38)(23,46)(3,69)(13,39)(26,52)(17,51)(34,68)(5,65)
    (35,70)(11,55)(33,66)(22,44)(21,42)(9,63)(27,54)(15,45)(25,75)(10,50)(14,56)(4,28)
    (20,60)(6,30)(8,40)(18,72)(12,36)(24,48)(16,64) 32

76  (29,58)(31,62)(37,74)(7,49)(23,69)(2,46)(13,26)(3,39)(17,51)(34,68)(4,52)(19,57)(38,76)
    (5,65)(35,70)(11,55)(33,66)(22,44)(21,42)(9,63)(27,54)(15,45)(25,75)(10,50)(14,28)
    (20,40)(6,36)(8,56)(18,72)(30,60)(12,24)(16,48)(32,64)

77  (29,58)(31,62)(37,74)(7,49)(11,77)(5,55)(35,70)(13,65)(3,39)(33,66)(17,51)(19,57)(23,69)
    (2,46)(22,44)(26,52)(34,68)(38,76)(21,42)(9,63)(27,54)(15,45)(25,75)(10,50)(14,56)(4,28)
    (20,60)(6,30)(8,40)(18,72)(12,36)(24,48)(16,64) 32

78  (29,58)(31,62)(37,74)(7,49)(11,77)(5,55)(35,70)(13,65)(23,69)(2,46)(22,66)(3,33)(39,78)
    (26,52)(4,44)(17,51)(34,68)(19,57)(38,76)(21,42)(9,63)(27,54)(15,45)(25,75)(10,50)
    (14,28)(20,40)(30,60)(8,56)(18,72)(6,36)(12,24)(16,48)(32,64)

80  (29,58)(31,62)(37,74)(7,49)(11,77)(5,55)(35,70)(13,65)(23,69)(2,46)(22,66)(3,33)(39,78)
    (26,52)(4,44)(17,51)(34,68)(19,57)(38,76)(21,42)(9,63)(27,54)(15,45)(25,75)(10,50)
    (14,56)(6,30)(18,72)(12,36)(20,60)(8,24)(40,80)(16,48)(32,64) 28

81  (29,58)(31,62)(37,74)(7,49)(11,77)(5,55)(35,70)(13,65)(23,69)(2,46)(22,66)(3,33)(39,78)
    (26,52)(4,44)(17,51)(34,68)(19,57)(38,76)(21,42)(9,63)(15,45)(25,75)(10,50)(27,81)
    (14,28)(6,30)(18,54)(12,36)(20,60)(8,56)(24,72)(16,48)(32,64)(40,80)

84  (29,58)(31,62)(37,74)(41,82)(7,49)(11,77)(5,55)(35,70)(13,65)(23,69)(2,46)(22,66)(3,33)
    (39,78)(26,52)(4,44)(17,51)(34,68)(19,57)(38,76)(25,75)(10,50)(9,45)(21,63)(27,81)
    (15,60)(6,30)(18,54)(20,80)(8,40)(32,64)(12,84)(36,72)(24,48)(14,42)(28,56) 16

85  (29,58)(31,62)(37,74)(41,82)(7,49)(11,77)(5,55)(35,70)(13,65)(17,85)(3,51)(33,66)(39,78)
    (19,57)(23,69)(2,46)(22,44)(26,52)(34,68)(38,76)(25,75)(10,50)(9,45)(21,63)(27,81)
    (15,30)(6,12)(18,54)(42,84)(14,56)(4,28)(20,60)(8,16)(32,64)(36,72)(24,48)(40,80)

92 (31,62)(37,74)(41,82)(43,86)(7,49)(11,77)(5,55)(35,70)(13,65)(17,85)(3,51)(33,66)(39,78)
   (19,57)(23,69)(29,87)(2,58)(26,52)(34,68)(38,76)(46,92)(22,88)(4,44)(25,75)(10,50)(9,63)
   (27,81)(21,84)(28,56)(14,42)(45,90)(15,60)(6,30)(18,54)(20,80)(8,40)(16,48)(12,36)(24,72)
   (32,64) 91

110 (37,74)(41,82)(43,86)(47,94)(53,106)(29,87)(2,58)(31,62)(3,93)(19,57)(38,76)(23,69)
    (46,92)(5,95)(13,65)(39,78)(17,85)(51,102)(34,68)(7,91)(49,98)(11,77)(55,110)(26,104)
    (4,52)(33,66)(9,99)(21,63)(27,81)(22,44)(35,70)(15,105)(45,90)(25,75)(10,50)(20,100)
    (6,18)(12,84)(8,88)(40,80)(32,64)(14,42)(28,56)(24,48)(54,108)(36,72)(16,96)(30,60)

116 (41,82)(43,86)(47,94)(53,106)(31,93)(2,62)(37,74)(3,111)(23,69)(46,92)(29,87)(58,116)
    (5,115)(13,65)(39,78)(17,85)(51,102)(34,68)(7,91)(49,98)(11,77)(55,110)(19,95)(57,114)
    (38,76)(26,104)(4,52)(33,66)(9,99)(21,63)(27,81)(22,88)(35,70)(15,105)(45,90)(25,75)
    (10,50)(20,100)(6,18)(30,60)(8,40)(16,80)(32,64)(14,42)(28,84)(56,112)(24,48)(54,108)
    (36,72)(12,96) 44

117 (41,82)(43,86)(47,94)(53,106)(31,93)(2,62)(37,74)(3,111)(23,69)(46,92)(29,87)(58,116)
    (5,115)(13,65)(17,85)(51,102)(34,68)(7,91)(49,98)(11,77)(55,110)(19,95)(57,114)(38,76)
    (33,66)(9,99)(21,63)(27,81)(39,117)(22,88)(4,44)(35,70)(15,105)(45,90)(25,75)(10,50)
    (20,100)(6,18)(30,60)(26,78)(52,104)(8,40)(16,80)(32,64)(14,42)(28,84)(56,112)(24,48)
    (54,108)(36,72)(12,96)
\end{verbatim}
\normalsize
\bibliographystyle{amsplain}
\bibliography{juniper}

\end{document}